  \documentclass[12pt,onecolumn, final]{IEEEtran}
\IEEEoverridecommandlockouts                              
 \usepackage{theorem}

\usepackage{xcolor}
\usepackage{enumerate}
 \usepackage{graphics} 
\usepackage{epsfig} 
\usepackage{url}
\usepackage{amsmath} 
 \usepackage{amssymb}  
\usepackage{verbatim} 
\usepackage{setspace} 

\newcommand{\bea}{\begin{eqnarray}}
\newcommand{\eea}{\end{eqnarray}}
\newcommand{\beas}{\begin{eqnarray*}}
\newcommand{\eeas}{\end{eqnarray*}}
\newcommand{\be }{\begin{equation}}
\newcommand{\ee }{\end{equation}}

\newcommand{\bi}{\begin{itemize}}
\newcommand{\ei}{\end{itemize}}

\newcommand{\dd}{\mathrm{d}}

\newcommand{\dist}{\operatorname{dist}}
\newcommand{\diag}{\operatorname{diag}}

\newcommand{\Int}{\operatorname{{\mathrm Int}}}

\newcommand{\tr}{\operatorname{tr}}

\newcommand {\R}{\mathbb R}

\newtheorem{Theorem}{Theorem}
\newtheorem{Definition}{Definition}
\newtheorem{Lemma}{Lemma}
\newtheorem{Proposition}{Proposition}
\newtheorem{Corollary}{Corollary}

 \newtheorem{Example}{Example}
\newtheorem{Remark}{Remark}

\newcommand{\sname}{}
\newcommand{\slabel}[1]{\debug{\fbox{\tiny \sname #1}}\label{\sname #1}}
\newcommand{\debug}[1]{}              
\newcommand{\FB}{\begin{figure}[t]\centering}
\newcommand{\FE}[2]{\caption{#2
\debug{\fbox{\sname #1}}} \slabel{#1} \end{figure}}
\newcommand{\tB}{\begin{table}[hbtp]\centering}
\newcommand{\tE}[2]{\caption{#2
\debug{\fbox{\sname #1}}}\slabel{#1} \end{table}}

 \newcommand{\sostfull}{contractive  after a small overshoot and short transient}
 \newcommand{\sostshort}{SOST}

 \newcommand{\stfull}{contractive after a  short transient}
\newcommand{\stshort}{ST}

\newcommand{\soshort}{SO}
\newcommand{\sofull}{contractive after a  small overshoot}

\newcommand{\wefull}{weakly expansive}
\newcommand{\sweshort}{WE}
\title{\Large  Contraction After Small  Transients}

\author{Michael Margaliot and Eduardo D. Sontag  and    Tamir Tuller
\thanks{An abridged version of this paper was presented at the IEEE CDC 2014~\cite{3gen_cdc_2014}.
 MM (michaelm@eng.tau.ac.il) is with the School of Electrical Engineering-Systems
and the Sagol School of Neuroscience, Tel
Aviv University, Israel 69978;
EDS (eduardo.sontag@gmail.com) is with the
   Dept. of Mathematics and the   Center  for Quantitative Biology, Rutgers University, Piscataway, NJ 08854, USA;
TT (tamirtul@post.tau.ac.il) is with the Dept. of Biomedical Engineering and the Sagol
School of Neuroscience,  Tel-Aviv University, Tel-Aviv 69978, Israel.
EDS's work is supported in part by grants
           NIH  1R01GM100473, AFOSR FA9550-14-1-0060, and ONR N00014-13-1-0074.
 The research of MM and TT is partly supported by a research grant from the Israeli Ministry of Science,
Technology and Space.
}}

\begin{document}
\maketitle

\thispagestyle{plain}
\pagestyle{plain}

\begin{abstract}

Contraction theory is a powerful tool for proving asymptotic
properties of nonlinear dynamical systems including
convergence to an attractor  and  entrainment to a periodic excitation.
\textcolor{black}{We consider  three      generalizations
of   contraction with respect to a norm
that  allow
   contraction to take place
     after
small transients in time and/or amplitude.}
These    generalized contractive systems~(GCSs) are useful  for  several reasons.
First, we show that there exist simple and checkable conditions
guaranteeing that a system is a  GCS,  and
       demonstrate their usefulness
    using several
     models  from systems biology.
Second, allowing small transients
does not destroy the important
asymptotic properties of contractive systems \textcolor{black}{like convergence to a unique equilibrium point, if it exists, and entrainment to a periodic excitation}.
Third,  in some cases as we change the parameters in a
 contractive system it becomes a~GCS just before it looses \textcolor{black}{contractivity with respect to a norm}. In this respect, generalized contractivity  is the
analogue of marginal stability in Lyapunov stability theory.
%
\end{abstract}

\begin{keywords} Differential analysis,
contraction, stability, entrainment, phase locking, systems biology.
 \end{keywords}

\section{Introduction}
Differential analysis is based on studying the time evolution of the distance between trajectories
emanating from different initial conditions.
 A  dynamical  system is called contractive
  if  any two   trajectories   converge to one other at an exponential rate. This implies many desirable
  properties including convergence to a unique  attractor (if it exists),
  and entrainment to   periodic excitations~\cite{LOHMILLER1998683,entrain2011,sontag_cotraction_tutorial}.
Contraction theory proved to be a  powerful tool for analyzing nonlinear dynamical systems,  with applications in
control theory~\cite{cont_mech},
observer design~\cite{observer_posi_2011},
synchronization of coupled oscillators~\cite{wang_slotine_2005}, and more.
Recent extensions include: the notion of partial contraction~\cite{partial_cont},
analyzing   networks of interacting
   agents using contraction theory~\cite{russo_hier,Arcak20111219},
 a Lyapunov-like characterization of incremental stability~\cite{angeli_inc},
 and
a LaSalle-type  principle for contractive systems~\cite{contra_sep}.
\textcolor{black}{There is also a growing interest
in design techniques providing controllers that  render control systems contractive or
incrementally
stable; see, e.g.~\cite{Zamani2013949} and the references therein,
 and also the incremental ISS condition in~\cite{Desoer_cont}). }

A contractive system with added
diffusion terms or random noise  still satisfies
certain asymptotic properties~\cite{Aminzare201331,slotine_cont_noise}.
In this respect,
contraction is  a robust property.

In this paper, we introduce three forms of generalized contractive systems~(GCSs).
  These
  are motivated by requiring \textcolor{black}{contraction with respect to a norm}
	to take place
   only after
  arbitrarily small
transients in time and/or amplitude. Indeed,   contraction is usually
used to prove \emph{asymptotic} properties, and thus allowing (arbitrarily small)
transients seems reasonable.
In some cases as we change the parameters in a
 contractive system it becomes a~GCS just before it looses contractivity. In this respect, a GCS is the
analogue of marginal stability in Lyapunov stability theory.
We provide several sufficient conditions for a system to be a~GCS.
These conditions are   checkable, and we
  demonstrate their usefulness  using
    examples
 of systems that are \emph{not} contractive with respect to any norm, yet are
 GCSs.


 We begin with a brief review of some
 ideas from contraction theory.
 For more details, including the historic development of contraction theory,
 and the relation to other notions, see e.g.~\cite{soderling_survey,cont_anc,RufferWouwMueller:2013:Convergent-Systems-vs.-Incremental-Stabi:}.

 Consider the time-varying
 system
 \be\label{eq:fdyn}
            \dot{x}=f(t,x),
 \ee
with the state $x$ evolving on a \textcolor{black}{positively invariant}
convex set~$\Omega \subseteq \R^n $. We assume that~$f(t,x)$ is differentiable with respect to~$x$, and that
both~$f(t,x)$ and~$J(t,x):=\frac{\partial f}{\partial x}(t,x)$ are continuous in~$(t,x)$.
Let~$x(t,t_0,x_0)$   denote the solution of~\eqref{eq:fdyn}
 at time~$t \geq t_0$ with~$x(t_0)=x_0$
 (for the sake of simplicity, we assume from here on
 that
 $x(t,t_0,x_0)$ exists and is unique for all~$t \geq t_0\geq 0$
 and all~$x_0 \in \Omega$).

We say that~\eqref{eq:fdyn} is
  \emph{contractive} on~$\Omega$ with respect to  a norm~$|\cdot| :\R^n \to \R_+$
 if there exists~$c>0$
 such that
 \be\label{eq:contdef}
            |x(t_2,t_1,a)-x(t_2,t_1,b)|  \leq   \exp( -   (t_2-t_1) c ) |a-b|
 \ee
 for all $t_2\geq t_1 \geq   0 $ and all~$a,b \in \Omega$.
 In other words,   any two trajectories contract to one another at an exponential rate. This implies
 in particular that the initial condition is ``quickly forgotten''. \textcolor{black}{Note that Ref.~\cite{LOHMILLER1998683}
provides a more general and intrinsic definition, where contraction is with respect to a time- and state-dependent metric~$M(t,x)$ (see also~\cite{SimpsonPorco201474} for a general treatment of contraction on a Riemannian
manifold).  Some of the results below may be stated  using this more general framework.
But, for a given dynamical system finding such a metric may be difficult.
Another extension  of  contraction is
incremental stability~\cite{angeli_inc}. Our approach is based on the
fact that there 	exists a  simple {sufficient} condition guaranteeing~\eqref{eq:contdef},
so generalizing~\eqref{eq:contdef} appropriately
leads to \emph{checkable} sufficient  conditions for a system to  be a  GCS.
Another advantage of our approach  is that  a GCS retains the important property of
  entrainment  to periodic signals. }

Recall that a vector norm~$|\cdot|:\R^n \to \R_+$
 induces a matrix measure
$\mu :\R^{n\times n} \to \R$ defined by
$
            \mu(A):=\lim_{\epsilon \downarrow 0} \frac{1}{\epsilon}
            (||I+\epsilon A||  -1),
$
where~$||\cdot||:\R^{n\times n}\to \R_+$ is the matrix norm induced by~$|\cdot|$.
A standard  approach for
proving~\eqref{eq:contdef} is based on  bounding some matrix measure of
 the Jacobian~$J$. Indeed, it is well-known (see, e.g.~\cite{entrain2011})
that if  there exist a vector norm~$|\cdot|$ and~$c>0$ such that the induced
 matrix measure~$\mu:\R^{n\times n} \to \R$ satisfies
 \be\label{eq:jtc}
            \mu(J(t,x)) \leq -c,
\ee
for all $t_2\geq t_1 \geq   0 $ and all~$x \in \Omega$
then~\eqref{eq:contdef}
holds. (This is in fact a particular case of using a Lyapunov-Finsler
 function to prove contraction~\cite{contra_sep}.)

 It is well-known~\cite[Ch.~3]{vid} that the matrix measure induced by the~$L_1$ vector norm is
\be\label{eq:muqdef}
 \mu_1(A)=\max\{c_1(A),\ldots, c_n(A)\},
\ee
 where
\be\label{eq:ccstac}
      c_j(A):=A_{jj}+\sum_{ \substack { 1\leq i \leq n\\  i \not = j} } |A_{ij}|  ,
 \ee
  i.e.,   the sum of the entries in column~$j$ of~$A$,
 with non diagonal
elements   replaced by their absolute values.
The
 matrix measure induced by the~$L_\infty$ norm is
\be\label{eq:mat_meas_inf}
 \mu_\infty(A)=\max\{d_1(A),\ldots, d_n(A)\},
\ee
 where
\be\label{eq:dstac}
      d_j(A):=A_{jj}+\sum_{ \substack { 1\leq i \leq n\\  i \not = j} } |A_{ji}|  ,
\ee
  i.e.,   the sum of the entries in row~$j$ of~$A$,
 with non diagonal
elements   replaced by their absolute values.

Often it is useful to work with scaled norms.
 Let~$|\cdot |_*$ be some  vector norm, and let~$\mu_*:\R^{n\times n}\to\R$
denote its induced matrix measure. If~$P\in\R^{n\times n}$ is an invertible matrix, and
$|\cdot|_{*,P} : \R^n \to \R_+$ is the vector norm defined by
$|z|_{*,P}:=|P  z|_* $ then the
induced matrix measure is $
                    \mu_{*,P}(A) = \mu_*(PAP^{-1}).$

One important  implication of
contraction is \emph{entrainment} to a periodic excitation.
Recall that~$f:\R_+ \times  \Omega  \to \R^n$ is called~\emph{$T$-periodic} if
 \[
 f(t,x) = f(t + T,x)
  \]
  for all $t \geq 0$ and all~$x \in \Omega$. \textcolor{black}{Note that for the system~$\dot x(t)=f(u(t),x(t))$,
	with~$u$ an input (or excitation) function,~$f$ will be~$T$ periodic if~$u$ is
	a~$T$-periodic function.} It is well-known \textcolor{black}{\cite{LOHMILLER1998683, entrain2011}}
	that
 if~\eqref{eq:fdyn} is contractive
and~$f$  is~$T$-periodic
then  for any~$t_1\geq 0$ there exists  a unique periodic solution $\alpha:[t_1,\infty) \to \Omega$ of~\eqref{eq:fdyn},
of period~$T$, and every trajectory  converges to~$\alpha$.
Entrainment   is important  in various applications
ranging from biological systems~\cite{RFM_entrain,entrain2011}
to the stability of  a  power grid~\cite{dorf-bullo}.
Note that for the particular case  where~$f$ is time-invariant,  this  implies that
if~$\Omega$ contains an equilibrium point~$e$
 then it is unique and all  trajectories converge to~$e$.

The remainder of this paper   is organized as follows.
\textcolor{black}{Section~\ref{sec:defe} presents three generalizations of~\eqref{eq:contdef}}. Section~\ref{sec:main}
 details sufficient conditions for their existence,
and describes their implications.    The proofs of all the  results are detailed
in Section~\ref{sec:proofs}.

\section{Definitions of contraction after small transients}\label{sec:defe}
We begin by defining  three generalizations of~\eqref{eq:contdef}.
 \begin{Definition}\label{def:qcont}
 The time-varying system~\eqref{eq:fdyn} is said to be:
\begin{itemize}
\item  \emph{\sostfull}~(\sostshort) on~$\Omega$
w.r.t.  a norm~$|\cdot| :\R^n \to \R_+$ if for each~$\varepsilon >0$
and each~$\tau>0$
  there exists~$\ell=\ell(\tau,\varepsilon)>0$
 such that
 \begin{align}\label{eq:qcont}
            |x(t_2+\tau,& t_1,a)-   x(t_2+\tau,t_1,b)|   \leq   (1+\varepsilon) \exp(-  (t_2-t_1) \ell ) |a-b| \,
 \end{align}
for all $t_2\geq t_1\geq 0$  and all $a,b \in \Omega$.
 \item   \emph{\sofull} ({\soshort}) on~$\Omega$
w.r.t.  a norm~$|\cdot| :\R^n \to \R_+$ if for  each~$\varepsilon >0$
  there exists~$\ell=\ell(\varepsilon)>0$
 such that
 \begin{align}\label{eq:uninew}
            |x(t_2, & t_1,a)-x(t_2,t_1,b)|
             \leq   (1+\varepsilon) \exp(-  (t_2-t_1) \ell ) |a-b|
 \end{align}
 for all  $t_2\geq t_1\geq 0$ and all  $a,b \in \Omega$.
\item    \emph{\stfull} (\stshort) on~$\Omega$
w.r.t. a norm~$|\cdot| :\R^n \to \R_+$ if for
 each~$\tau>0$
  there exists~$\ell=\ell(\tau)>0$
 such that
 \begin{align}\label{eq:ucont}
            |x(t_2+\tau, & t_1,a)-x(t_2+\tau,t_1,b)|
            \leq     \exp(-  (t_2-t_1) \ell ) |a-b|
 \end{align}
  for all $t_2\geq t_1\geq 0$  and all $a,b \in \Omega$.
\end{itemize}
 \end{Definition}

The definition of {\sostshort}  is motivated
by  requiring   contraction  at an exponential rate,
    but only
    after an (arbitrarily small) time~$\tau$, and with  an (arbitrarily small)   overshoot~$(1+\varepsilon)$.
    However, as we will see below when the convergence rate~$ \ell$
    may depend  on~$\varepsilon$
      a somewhat  richer  behavior may occur.
The definition of~{\soshort} is similar to that of~\sostshort, yet now the convergence
rate~$\ell$ depends only on~$\varepsilon$, and there is
no time transient~$\tau$ (i.e.,~$\tau=0$). In other words,~{\soshort} is a uniform (in~$\tau$)
version of~\sostshort.
The third definition,~{\stshort},  allows the contraction to ``kick in''
only after a time transient of length~$\tau$.

It is clear    that every contractive
 system is~\sostshort,
{\soshort}, and~{\stshort}. Thus, all these notions
  are
    generalizations of   contraction.
Also, both  {\soshort} and {\stshort}   imply~{\sostshort} and, as
we will see
below, under a mild technical
condition on~\eqref{eq:fdyn} {\soshort}
and   {\sostshort}  are equivalent.
Figure~\ref{fig:graphbn2} on p.~\pageref{fig:graphbn2}
  summarizes  the relations between these  GCSs (as well as  other notions defined below).

  The motivation for these definitions stems from the fact that   important
    applications of contraction are  in proving \emph{asymptotic} properties.
    For example,
    proving that an equilibrium point is globally attracting
    or that the state-variables entrain to a  periodic excitation.
    These properties describe what happens as~$t\to \infty$, and so
    it seems natural to generalize contraction in a way that allows
    initial  transients in time and/or amplitude.

The next  simple  example demonstrates a system that does not satisfy~\eqref{eq:contdef}, but is a GCS.
\begin{Example}\label{exa:scalarsys}
            Consider the \emph{scalar} time-varying  system
            \be\label{eq:scals}
                    \dot x(t)=-\alpha(t)x(t),
            \ee
             with the state $x$ evolving   on~$\Omega:=[-1,1]$, and~$\alpha:\R_+\to\R_+$ is a class~K function (i.e.~$\alpha$ is continuous and strictly increasing, with~$\alpha(0)=0$).
           It is straightforward to
           show that this system does not satisfy~\eqref{eq:contdef}
           w.r.t. \emph{any} norm
           (note that the Jacobian~$J(t)=-\alpha(t)$  satisfies~$J(0)=0$),
					yet it is~{\stshort}, with $\ell(\tau) = \alpha(\tau)>0$, for any given $\tau>0$.
\end{Example}

The next section presents our main results. The proofs are placed in Section~\ref{sec:proofs}.

\section{Main Results}\label{sec:main}

The next three subsections study the three forms of~GCSs defined above.
\subsection{{\sostfull}}

Just like contraction,~{\sostshort} implies entrainment to a
periodic excitation. To show this,
assume that the vector field~$f$ in~\eqref{eq:fdyn}
	is~$T$ periodic.
Pick~$t_0\geq 0$.
 Define~$m:\Omega \to \Omega$ by~$m(a):=x(T+t_0,t_0,a)$.
 In other words,~$m$ maps~$a$
to the solution of~\eqref{eq:fdyn} at time~$T+t_0$ for the initial condition~$x(t_0)=a$.
Then~$m$ is continuous and  maps the convex and compact set~$\Omega$ to itself,
 so
by the Brouwer fixed point theorem (see, e.g.~\cite[Ch.~6]{fixedpoint}) there exists~$\zeta \in \Omega$ such that~$m(\zeta)=\zeta$, i.e.~$x(T+t_0,t_0,\zeta)=\zeta$.
 This implies that~\eqref{eq:fdyn}
 admits a periodic solution~$\gamma:[t_0,\infty) \to \Omega$
  with period~$T$. Assuming that the system is also~{\sostshort}, pick~$\tau,\varepsilon>0$.
	Then there exists~$\ell=\ell(\tau,\varepsilon)>0$   such that
\[
			|x(t-t_0+\tau,t_0,a)- x(t   - t_0+\tau,t_0,\zeta)|\leq(1+\varepsilon)\exp(-(t - t_0) \ell)|a-\zeta|,
\]
for all~$a\in \Omega$ and all~$t\geq t_0$. Taking~$t\to\infty$ implies that every solution converges to~$\gamma$.
 In particular,   there cannot be two distinct periodic solutions.
Thus, we proved the following.

\begin{Proposition}\label{prop:st_entrain}
 Suppose that the time-varying system~\eqref{eq:fdyn}, with state~$x$
 evolving on a compact and convex state-space~$\Omega\subset\R^n$,
is {\sostshort},
and that the vector field~$f$ is~$T$-periodic.
Then for any~$t_0 \geq 0$ it admits
a unique periodic solution~$\gamma:[t_0,\infty ) \to\Omega$  with period~$T$, and~$x(t,t_0,a)$ converges to~$\gamma$ for any~$a\in\Omega$.
\end{Proposition}

Since both~{\soshort} and~{\stshort} imply~{\sostshort}, Proposition~\ref{prop:st_entrain} holds for all three forms of~GCSs.

\textcolor{black}{Our next goal is to derive a sufficient condition for~{\sostshort}.
One may naturally expect that if~\eqref{eq:fdyn} is contractive w.r.t. a set of norms~$|\cdot|_\zeta$,
with, say~$\zeta \in (0,p]$, $p>0$, and  that~$\lim_{\zeta\to 0}|\cdot|_\zeta=|\cdot|$ then
\eqref{eq:fdyn}  is a GCS w.r.t.  the norm~$|\cdot|$. In fact, this can be further generalized by requiring
\eqref{eq:fdyn} to be contractive w.r.t.~$|\cdot|_\zeta$ only on suitable subset~$\Omega_\zeta$ of the state-space.
This leads to  the following definition.
}

 \begin{Definition}\label{def:NC}
 System~\eqref{eq:fdyn} is said
  to be \emph{nested contractive}~(NC) on~$\Omega$ with respect to a norm~$|\cdot|$
if there exist  convex   sets~$\Omega_\zeta \subseteq \Omega$,    and norms~$|\cdot|_\zeta:\R^n\to \R_+$,
 where~$\zeta\in (0,1/2]$,
 such that the following conditions hold.
 \begin{enumerate}[(a)]
                    \item $\cup_{\zeta \in (0,1/2]} \Omega_\zeta=\Omega$, and
                    \be\label{eq:setsinc}
                    \Omega_{\zeta_1} \subseteq \Omega_{\zeta_2},\quad \text{for all } \zeta_1 \geq \zeta_2.
                    \ee
                    \item For every~$ \tau>0$ there exists~$\zeta=\zeta( \tau)\in(0,1/2]$,
                    with~$\zeta ( \tau)\to 0$ as~$ \tau \to 0 $,
                    such that  for every~$a\in \Omega$ and every~$t_1 \geq 0$
                    \be \label{eq:enter}
                            x(t ,t_1,a)\in\Omega_{\zeta },\quad \text{for all } t\geq t_1+\tau,
                     \ee
                     and~\eqref{eq:fdyn}
                                              is contractive  on~$\Omega_\zeta$ with respect to~$|\cdot|_\zeta$.
                     \item \label{item:cc} The  norms  $|\cdot|_\zeta$ converge to~$|\cdot| $ as~$\zeta\to 0$, i.e., for every~$\zeta>0$ there exists~$s =s(\zeta) >0$, with~$s(\zeta) \to 0$
                          as~$ \zeta \to 0$, such that
                         \[
                                  (1-s)|  y | \leq |y|_\zeta  \leq (1+s) |y|      ,\quad \text{for all } y \in \Omega       .
                         \]
 \end{enumerate}
 \end{Definition}

Eq.~\eqref{eq:enter} means that
 after an arbitrarily  short time every trajectory enters and remains in
 a subset~$\Omega_\zeta$ of the state space on which we have contraction
 with respect to~$|\cdot|_\zeta$. \textcolor{black}{We can now state the main result in this subsection.}

 \begin{Theorem}\label{thm:qcon}
 If the system~\eqref{eq:fdyn} is NC  w.r.t. the norm~$|\cdot|$  then it is
  {\sostshort} w.r.t. the norm~$|\cdot|$.
 \end{Theorem}

The next example  demonstrates Theorem~\ref{thm:qcon}. It also shows that
as we change the parameters in a contractive system, it may become a~GCS
when it hits the ``verge'' of contraction \textcolor{black}{(as defined in~\eqref{eq:contdef})}.
This is reminiscent of an
asymptotically stable  system that becomes
marginally stable  as it looses stability.

\begin{Example}\label{exa:bio_smith}
Consider the system
\begin{align}\label{eq:bio}
\dot{x}_1&=g(x_n)-\alpha_1 x_1,\nonumber\\
\dot{x}_2&=x_1-\alpha_2 x_2\nonumber,\\
\dot{x}_3&=x_2-\alpha_3 x_3 \nonumber,\\
       &\vdots\nonumber\\
\dot{x}_n&=x_{n-1}-\alpha_n x_n  ,
\end{align}
where~$\alpha_i>0$, and
\[
g(u):=\frac{1+u}{k+u} , \quad \text{with } k>1.
\]
As explained  in~\cite[Ch.~4]{hlsmith} this may model a simple
  biochemical     feedback
 control circuit for protein synthesis in the cell.
The~$x_i$s represent concentrations of various macro-molecules in the cell
and therefore must be non-negative. It is straightforward to
verify that~$x(0) \in \R^n_+$ implies that~$x(t) \in \R^n_+$ for all~$t \geq 0$.

Let~$\alpha:=\prod_{i=1}^n \alpha_i  $, and for~$\varepsilon>0$ let
\begin{align*}
D_\varepsilon :=  \diag \left(1,\alpha_1-\varepsilon , (\alpha_1-\varepsilon)(\alpha_2-\varepsilon ),\dots, \prod_{i=1}^{n-1}(\alpha_{i}-\varepsilon ) \right).
\end{align*}
We show in Section~\ref{sec:proofs} that if
 \be\label{eq:km1}
   {k-1}  < \alpha {k^2}
   \ee
then~\eqref{eq:bio} is contractive  on~$\R_+^n$ w.r.t.  the scaled
norm~$|\cdot|_{1,D_\varepsilon}$ for all~$\varepsilon>0$ sufficiently small.
If
$ {k-1} =\alpha {k^2}
 $  then~\eqref{eq:bio} does not satisfy~\eqref{eq:contdef}, w.r.t.  any norm,
on~$\R_+^n$, yet it is~{\sostshort} on~$\R_+^n$ w.r.t.
 the norm~$|\cdot|_{1,D_0}$.

Note that
 for all~$x\in\R^n_+$,
\be\label{eq:gderi}
        g'(x_n)= \frac{k-1}{(k+x_n)^2} \leq \frac{k-1}{k^2} =g'(0).
\ee
Thus~\eqref{eq:km1} implies that the system  satisfies~\eqref{eq:contdef}
if and only
if the
``total dissipation'' $\alpha$
is strictly larger   than~$g'(0)$.

Using the fact that~$g(u)<1$ for all~$u\geq 0$ it is straightforward to show that
the set
\[
        \Omega_r:= r  ( [0,\alpha_1^{-1}]\times[0, (\alpha_1 \alpha_2)^{-1}] \times\dots\times [0, \alpha^{-1}]  )
\]
is an invariant set of the dynamics for all~$r\geq 1$.
Thus,~\eqref{eq:bio},
with~$   {k-1}  \leq \alpha {k^2}$,
admits a unique equilibrium point~$e \in \Omega_1 $ and
\[
            \lim_{t\to \infty} x(t,a)=e,\quad \text{for all } a \in \R^n_+.
\]
This property also follows from a more general result~\cite[Prop.~4.2.1]{hlsmith}
that is proved using the theory of irreducible  cooperative  dynamical systems.
Yet the contraction approach leads to new insights.
For example,
it implies that the distance between trajectories can only decrease,
and can also be used to prove entrainment to   suitable generalizations
of~\eqref{eq:bio} that include periodically-varying inputs.
\end{Example}

\textcolor{black}{
Cells often respond to external stimulus by
 modification of proteins. One   mechanism  for this
 is \emph{phosphorelay}  (also called phosphotransfer)
 in which a phosphate group is transferred through a serial
 chain of proteins
  from an initial  histidine
kinase~(HK)    down to a
final response regulator~(RR).
  The next example uses Theorem~\ref{thm:qcon}
 to analyze
 a model for phosphorelay from~\cite{phos_relays}.
}

\begin{Example}\label{exa:phos:relay}
Consider the system
\begin{align}\label{eq:phos}
								\dot x_1&= (p_1-x_1)c-\eta_1x_1(p_2-x_2),\nonumber\\
          			\dot x_2&=\eta_1x_1(p_2-x_2)-\eta_2x_2(p_3-x_3),\nonumber\\
								&\vdots\nonumber\\
								\dot x_{n-1}&=\eta_{n-2} x_{n-2}(p_{n-1}-x_{n-1})-\eta_{n-1} x_{n-1}(p_n-x_n),\nonumber\\
								\dot x_n&= \eta_{n-1} x_{n-1}(p_n-x_n) -\eta_n x_n,
\end{align}
where~$\eta_i,p_i>0$, and~$c:[t_1,\infty)\to \R_+$.
In the context of phosphorelay~\cite{phos_relays},~$c(t)$ is the strength at time~$t$ of the
stimulus activating the  HK,
$x_i(t)$ is the concentration of
the phosphorylated form of the protein at the $i$th layer at time~$t$, and
$p_i$ denotes the  total protein concentration at that layer. Note that
$\eta_n x_n$
is the  flow of
 the phosphate group to an external  receptor molecule.

In  the particular case where~$p_i=1$ for all~$i$~\eqref{eq:phos}
  becomes the \emph{ribosome flow model}~(RFM) \cite{reuveni}. This  is the mean-field approximation of an important model
	from non-equilibrium statistical physics called the \emph{totally asymmetric simple  exclusion process}~(TASEP) \cite{solvers_guide}.
  In the RFM,
 $x_i \in[0,1]$ is the normalized occupancy at  site~$i$, where~$x_i=0$ [$x_i=1$]
means that site~$i$ is completely free [full], and~$\eta_i$ is the capacity of the link that
connects site~$i$ to site~$i+1$. This has been used to model mRNA translation, where every site corresponds to a group of
codons on the mRNA strand,~$x_i(t)$
is the normalized occupancy   of ribosomes at site~$i$ at time~$t$,
$c(t)$ is the initiation rate at time~$t$, and $\eta_i$
 is the   elongation rate from site~$i$ to site~$i+1$.

Our original motivation for   generalizing~\eqref{eq:contdef}
        was to
prove entrainment in the  RFM~\cite{RFM_entrain}.
For more results on the~RFM, see~\cite{RFM_stability,RFM_feedback,HRFM_steady_state,infi_HRFM,rfm_concave}.

Assume that there exists~$\eta_0>0$ such that~$c(t)\geq\eta_0$ for all~$t\geq t_1$.
Let~$\Omega:=[0,p_1]\times\dots\times[0,p_n]$ denote the state-space of~\eqref{eq:phos}.
Then, as shown in Section~\ref{sec:proofs},~\eqref{eq:phos}  does not satisfy~\eqref{eq:contdef}, w.r.t. any norm,
on~$\Omega$, yet it is~{\sostshort} on~$\Omega$ w.r.t.   the~$L_1$ norm.

\end{Example}

Considering Theorem~\ref{thm:qcon} in the special case where all the sets~$\Omega_\zeta$ in Definition~\ref{def:NC}
are equal to~$\Omega$ yields the following result.
\begin{Corollary}\label{coro:new}
Suppose that  \eqref{eq:fdyn}  is contractive  on~$\Omega$ w.r.t. a
set of norms~$|\cdot|_\zeta$,~$\zeta\in (0,1/2]$,
                    and that
										condition~(\ref{item:cc}) in Definition~\ref{def:NC} holds.
Then~\eqref{eq:fdyn} is {\sostshort} on~$\Omega$ w.r.t.~$|\cdot|$.
 \end{Corollary}

Corollary~\ref{coro:new} may be  useful in cases where some matrix measure
of the Jacobian~$J$ of \eqref{eq:fdyn} turns out to be non positive on~$\Omega$,
 but not strictly negative, suggesting that the system is
``on the verge'' of satisfying~\eqref{eq:contdef}.
The next result demonstrates this for the
time-invariant system
\be\label{eq:time_in_var_sys}
                    \dot{x}=f(x),
\ee
and the
particular
case of the matrix measure~$\mu_1:\R^{n\times n}\to\R$ induced by the~$L_1$ norm.
Recall that this is given by~\eqref{eq:muqdef} with the~$c_j $s defined in~\eqref{eq:ccstac}.

\begin{Proposition}\label{prop:new_meas_zero}
Consider  the Jacobian~$J(\cdot):\Omega \to \R^{n\times n}$ of the time-invariant system~\eqref{eq:time_in_var_sys}.
Suppose that~$\Omega$ is compact and that the set~$\{1,\dots,n\}$
can be divided into two non-empty disjoint sets~$S_0$ and~$S_{-}$ such that the following properties hold
for all~$x \in \Omega$:
\begin{enumerate}
\item
for any~$k\in S_0$, $c_k(J(x))\leq 0$; \label{item:s0}
\item
for any~$j \in S_-$, $c_j(J(x))< 0$; \label{item:sminus}
\item
for any~$i \in S_0$ there exists an index~$z=z(i)\in S_-$ such that~$J_{zi}(x)>0$.
 \label{item:rec}
\end{enumerate}
Then~\eqref{eq:time_in_var_sys}  is {\sostshort}  on~$\Omega$ w.r.t. the~$L_1$ norm.
 \end{Proposition}

The proof of  Proposition~\ref{prop:new_meas_zero} is based on the following idea.
By compactness of~$\Omega$,
there exists~$\delta>0$ such that
\be \label{eq:comp_ass}
c_j(J(x))<-\delta,\quad \text{for all } j \in S_- \text{ and all }x\in\Omega.
\ee
The conditions stated  in the proposition imply
that there exists a diagonal matrix~$P$ such that~$c_k(PJ P^{-1})<0$
 for all~$k \in S_0$. Furthermore, there exists such a~$P$    with diagonal entries \emph{arbitrarily close}
 to~$1$,  so~$c_j(PJP^{-1})<-\delta/2$
 for all~$j \in S_-$. Thus,~$\mu_1(PJP^{-1})<0$. Now Corollary~\ref{coro:new} implies {\sostshort}.
\textcolor{black}{ Note that this implies that  the compactness assumption may be dropped
 if for example it is known that~\eqref{eq:comp_ass} holds.
}

\begin{Example}\label{exa:new_from_russo}
Consider the system:
\begin{align} \label{eq:trans_module}
\dot x=& -\delta x+k_1 y-k_2(e_T-y)x, \nonumber \\
\dot y=& - k_1 y +k_2(e_T-y)x,
\end{align}
where~$\delta,k_1,k_2,e_T>0$, and~$\Omega:=[0,\infty)\times[0,e_T]$.
This is a basic model for a transcriptional module that is  ubiquitous in both
biology and synthetic biology
(see, e.g.,~\cite{retro_2008,entrain2011}).
Here~$x(t)$ is the concentration at time~$t$ of a transcriptional factor~$X$
that regulates  a downstream
transcriptional module by binding to a promoter
 with concentration~$e(t)$ yielding a protein-promoter complex~$Y$ with concentration~$y(t)$.
 The binding reaction   is reversible with  binding and
dissociation rates~$ k_2$ and $k_1$, respectively. The linear degradation rate of~$X$ is~$\delta$,
and as the promoter is not
subject to decay, its total concentration, $e_T$, is conserved,
so~$e(t)=e_T-y(t)$.
The  Jacobian of~\eqref{eq:trans_module}
is~$J=\begin{bmatrix}    -\delta-k_2(e_T-y)& k_1+k_2 x \\
k_2(e_T-y) &-k_1-k_2 x\end{bmatrix}$,
and all the properties in Prop.~\ref{prop:new_meas_zero}
 hold with~$S_-=\{1\}$ and~$S_0=\{2\}$.
  Indeed,~$J_{12}=k_1+k_2x>k_1>0$ for all~$\begin{bmatrix} x& y \end{bmatrix}^T \in\Omega$.
Thus, \eqref{eq:trans_module} is
 {\sostshort} on~$\Omega$ w.r.t. the~$L_1$ norm.
\textcolor{black}{Note that Ref.~\cite{entrain2011}  showed that~\eqref{eq:trans_module} is
contractive w.r.t.  a certain \emph{weighted}  $L_1$ norm.
Here we showed  {\sostshort} w.r.t. the (unweighted) $L_1$ norm.}
\end{Example}

\begin{Example}\label{exa:sec_new_from_russo}
A more general example studied in~\cite{entrain2011} is where the transcription
factor regulates several independent downstream transcriptional modules.
This leads to the following model:
\begin{align} \label{eq:mult_trans_module}
\dot x=& -\delta x+k_{11} y_1-k_{21}(e_{T,1}-y_1)x+
k_{12} y_2-k_{22}(e_{T,2}-y_2)x+\dots+
k_{1n} y_n-k_{2n}(e_{T,n}-y_n)x, \nonumber \\
\dot y_1=& - k_{11} y_1 +k_{21}(e_{T,1}-y_1)x,\nonumber\\
&\vdots\nonumber \\
\dot y_n=& - k_{1n} y_n +k_{2n}(e_{T,n}-y_n)x,
\end{align}
where~$n$ is the number of regulated modules.
The state-space is~$\Omega=[0,\infty)\times[0,e_{T,1}]\times\dots\times [0,e_{T,n}]$.
 The  Jacobian of~\eqref{eq:mult_trans_module}
is
\[
J=\begin{bmatrix}    -\delta-\sum_{i=1}^n k_{2i}(e_{T,i}-y_i) & k_{11}+k_{21} x & k_{12}+k_{22} x & \dots & k_{1n-1}+k_{2n-1}x & k_{1n}+k_{2n}x  \\
k_{21}(e_{T,1}-y_1) &-k_{11}-k_{21} x & 0 &\dots& 0&0 \\
k_{22}(e_{T,2}-y_2) &0                &-k_{12}-k_{22} x & 0 &\dots &0\\
&\vdots\\
k_{2n}(e_{T,n}-y_n) &0                &0 & \dots&0 &-k_{1n}-k_{2n} x\\
\end{bmatrix},
\]
and all the properties in Prop.~\ref{prop:new_meas_zero}
 hold with~$S_-=\{1\}$ and~$S_0=\{2,3,\dots,n\}$.
Thus, this system is
 {\sostshort} on~$\Omega$ w.r.t. the~$L_1$ norm.
\end{Example}

 Arguing as in the proof of Proposition~\ref{prop:new_meas_zero} for the matrix measure~$\mu_\infty$
induced by the~$L_\infty$ norm (see~\eqref{eq:dstac})
yields the following result.

\begin{Proposition}\label{prop:new_meas_zero_infty}
Consider  the Jacobian~$J(\cdot):\Omega \to \R^{n\times n}$ of the time-invariant system~\eqref{eq:time_in_var_sys}.
 Suppose that~$\Omega$ is compact  and that  the set~$\{1,\dots,n\}$
can be divided into two non-empty disjoint sets~$S_0$ and~$S_{-}$ such that the following properties hold for all~$x\in\Omega$:
\begin{enumerate}
\item
$d_j(J(x))\leq0$ for all~$j \in S_0$; \label{item:s0_inf}
\item
$d_k(J(x))<0$ for all~$k \in S_{-}$; \label{item:sminus_inf}
\item
for any~$j \in S_0$
there exists an index~$z=z(j) \in S_{-}$ such that
$J_{jz}(x)\not = 0$. \label{item:rec_inf}
\end{enumerate}
Then~\eqref{eq:time_in_var_sys}  is {\sostshort}  on~$\Omega$ w.r.t. the~$L_\infty$ norm.
 \end{Proposition}

\subsection{\sofull}
A natural question is under what conditions~{\soshort} and~{\sostshort} are equivalent. To address this issue,
we introduce the following definition.

\begin{Definition}\label{eq:defntr}
                            We say that~\eqref{eq:fdyn} is \emph{\wefull} (\sweshort)
                            if for each~$\delta>0$ there exists~$\tau_0>0$ such that for all~$a,b \in \Omega$ and all~$t_0\geq 0$
                            \be\label{eq:sep}
                                        | x(t ,t_0,a)-x(t ,t_0,b) |\leq(1+\delta)|a-b|,\quad\text{for all }t \in [t_0,t_0+\tau_0].
                            \ee
 \end{Definition}

 \begin{Proposition}\label{prop:sepimp}
            Suppose that \eqref{eq:fdyn} is   {\sweshort}.
            Then \eqref{eq:fdyn} is {\sostshort}  if and only if it is {\soshort}.
\end{Proposition}

\begin{Remark}\label{rem:Global_Lip}
Suppose that~$f$ in~\eqref{eq:fdyn} is Lipschitz globally in~$\Omega$
 uniformly in~$t$, i.e. there exists~$L>0$ such that
\[
               |  f(t,x)-f(t,y)|\leq L|x-y|, \quad \text{for all }x,y \in \Omega, \; t\geq 0.
\]
Then by Gronwall's Lemma (see, e.g.~\cite[Appendix~C]{sontag_textbook})
 \[
        | x(t ,t_0,a)-x(t ,t_0,b) |\leq \exp \left (L (t-t_0) \right )|a-b|,
 \]
  for all $t \geq t_0\geq 0$, and this implies that~\eqref{eq:sep} holds for~$\tau_0:=\frac{1}{L} \ln(1+\delta)>0$.
In particular, if~$\Omega$ is compact and~$f$ is periodic in~$t$ then~{\sweshort} holds under rather weak continuity arguments on~$f$.
\end{Remark}

\subsection{{\stfull}}

For  \emph{time-invariant}  systems whose state evolves
on  a convex and  compact set
 it is possible to give a simple  sufficient condition for~{\stshort}.
Let~$\Int(S)$ [$\partial S$] denote the interior [boundary] of a set~$S$.
We require the following definitions.
\begin{Definition}
We say that~\eqref{eq:fdyn} is
\emph{non expansive}~(NE) w.r.t.   a  norm~$|\cdot|$ if
for all~$a,b\in\Omega$ and all~$s_2 > s_1\geq 0$
   \be\label{eq:exp}
        |x(s_2,s_1,a)-x(s_2,s_1,b)|   \leq    |a-b|.
\ee
We say that~\eqref{eq:fdyn} is \emph{weakly contractive}~(WC) if~\eqref{eq:exp}
holds with~$\leq $ replaced by~$<$.
\end{Definition}

\begin{Definition}\label{def:ic}
The  time-invariant system~\eqref{eq:time_in_var_sys}
  with the state $x$ evolving on a     compact and convex set~$\Omega \subset \R^n$,
 is said to be \emph{interior contractive}~(IC) w.r.t. a norm~$|\cdot|:\R^n\to\R_+$
 if the
following properties hold:
\begin{enumerate}[(a)]
\item
for every~$x_0 \in \partial \Omega$,
\be\label{eq:cond_a_enu}
     x(t,x_0) \not \in  \partial \Omega, \quad \text{for all } t>0;
\ee
\item  for every~$x \in \Int(\Omega)$,
\be \label{eq:mucom_inv}
\mu(J(x))<0,
\ee
where~$\mu:\R^{n\times n}\to \R$ is the matrix measure induced by~$|\cdot|$.
\end{enumerate}
\end{Definition}

In other words, the  matrix measure   is negative   in the interior of~$\Omega$,
and the boundary of~$\Omega$ is ``repelling''.
Note that these conditions   do not necessarily imply
that the system satisfies~\eqref{eq:contdef}
  on~$\Omega$, as it is possible that~$\mu(J(x))=0$   for some~$x\in \partial \Omega$.
Yet,~\eqref{eq:mucom_inv} does   imply that~\eqref{eq:time_in_var_sys} is~NE on~$\Omega$. \textcolor{black}{We can now state the main result in this subsection.}

\begin{Theorem}\label{thm:time_invar}
If the system~\eqref{eq:time_in_var_sys}  is IC w.r.t. a norm~$|\cdot|$ then
 it is~{\stshort} w.r.t.~$|\cdot|$.
\end{Theorem}

The proof of this result is based on showing that IC implies that
  for each $\tau > 0$ there exists~$d=d(\tau) > 0$
 such  that
 \[
 \dist(x(t,x_0),\partial \Omega) \geq d ,\quad \text{for all } x_0 \in \Omega  \text{ and all } t \geq \tau,
 \]
and then using this to conclude that for any~$t\geq \tau$
all the trajectories of the system are  contained in a convex and compact set~$D\subset \Int(\Omega)$. In this set the system is contractive with rate~$c:=\max_{x\in D}\mu(J(x))<0$.
The next example, that is a variation of a system studied in~\cite{entrain2011}, demonstrates this reasoning.

\begin{Example}\label{exa:new_from_russo_inf}
Consider a transcriptional factor~$X$
that regulates  a downstream
transcriptional module by irreversibly binding, at a rate~$k_2>0$,  to a promoter~$E$
  yielding a protein-promoter complex~$Y$.
The promoter is not
subject to decay, so its total concentration, denoted by~$e_T>0$, is conserved.
  Assume also that~$X$ is obtained
from an inactive form~$X_0$, for example through a phosphorylation
reaction that  is catalyzed by a kinase with abundance~$u(t)$ satisfying~$u(t)\geq u_0>0$ for all~$t\geq 0$.
 The sum  of the  concentrations of
$X_0$, $X$, and~$Y$ is constant, denoted by~$z_T$, with~$z_T>e_T$.
Letting~$x_1(t),x_2(t)$ denote the concentrations of~$X,Y$    at time~$t$
yields the model
\begin{align} \label{eq:trans_module_inf}
\dot x_1=& (z_T-x_1-x_2) u   -\delta x_1 -k_2(e_T-x_2)x_1, \nonumber \\
\dot x_2=&     k_2(e_T-x_2)x_1,
\end{align}
with the state  evolving on~$\Omega:=[0,z_T]\times[0,e_T]$.
\textcolor{black}{
Here~$\delta\geq 0$ is the dephosphorylation  rate $X \to  X_0$.}
Let~$P:=\begin{bmatrix}1&1\\0&1 \end{bmatrix}$, and consider the matrix measure~$\mu_{ \infty,P}$.
A calculation yields
\begin{align*}
									\tilde J&:=PJP^{-1}\\&=  \begin{bmatrix}
-u   -\delta & \delta  \\
k_2(e_T-x_2)  & k_2(x_2 -x_1-e_T)
\end{bmatrix},
\end{align*}
so~$d_1(\tilde J)=-u  -\delta +|\delta | \leq -u_0<0$,
and
\begin{align*}
d_2(\tilde J)&=  k_2(x_2 -x_1-e_T)+|k_2(e_T-x_2)  |\\
&=-k_2 x_1.
\end{align*}
Letting~$S:= \{0\}\times[0,e_T]$, we conclude that~$\mu_{\infty,P}(x)<0$ for all~$ x  \in (\Omega \setminus S)$.
For any~$x\in S$,
$\dot x_1 = (z_T -x_2) u  \geq(z_T -e_T) u_0>0  $, and arguing as in the proof of Theorem~\ref{thm:time_invar}
(see Section~\ref{sec:proofs}),
we conclude that for any~$\tau>0$ there exists~$d=d(\tau)>0$ such that
\[
					x_1(t,a)\geq d,\quad \text{for all } a\in \Omega \text{ and all } t\geq\tau.
\]
In other words, after time~$\tau$ all the trajectories are contained in the
closed and convex set~$D=D(\tau):=[d,z_T]\times[0,e_T]$. Letting~$c:=c(\tau)=\max_{x\in D}\mu_{\infty,P}(J(x))$ yields~$c<0$ and
\[
			 |x(t+\tau,a)-x(t+\tau,b)|_{\infty,P}\leq\exp(ct)|a-b|_{\infty,P},\quad\text{for all }a,b\in\Omega \text{ and all }t>0,
\]
  so~\eqref{eq:trans_module_inf}
	is~{\stshort} w.r.t.~$|\cdot|_{\infty,P}$.
\end{Example}

As noted above, the introduction of  GCSs
  is motivated by the idea that contraction is used to prove asymptotic results,
so allowing
initial transients should increase the class of systems that can be analyzed while still
allowing to prove asymptotic results.
The next result demonstrates this.
\begin{Corollary}\label{coro:attract}
      If~\eqref{eq:time_in_var_sys} is~IC with respect to some norm
  then it admits
  a unique equilibrium point~$e \in  \Int( \Omega )$,
and~$\lim_{t \to \infty}x(t,a)= e $ for all~$a\in \Omega$.
\end{Corollary}

\begin{Remark}
The proof of Corollary~\ref{coro:attract}, given in the Appendix, is based on
Theorem~\ref{thm:time_invar}.
 Consider the
  \emph{variational system} (see, e.g.,~\cite{contra_sep}) associated with~\eqref{eq:time_in_var_sys}:
\begin{align}\label{eq:var_sys}
\dot x&=f(x),\nonumber\\
\dot {\delta x}&=J(x)\delta x.
\end{align}
 An alternative  proof of Corollary~\ref{coro:attract}
is possible, using the  Lyapunov-Finsler function~$V(x,\delta x):=|\delta x| $,
where~$|\cdot|:\R^n\to\R_+$ is the
vector norm corresponding to the matrix measure~$\mu$ in~\eqref{eq:mucom_inv},
 and the
 LaSalle invariance principle
  described  in~\cite{contra_sep}.
\end{Remark}

\begin{figure}[t]
  \begin{center}
  \includegraphics[height=7cm]{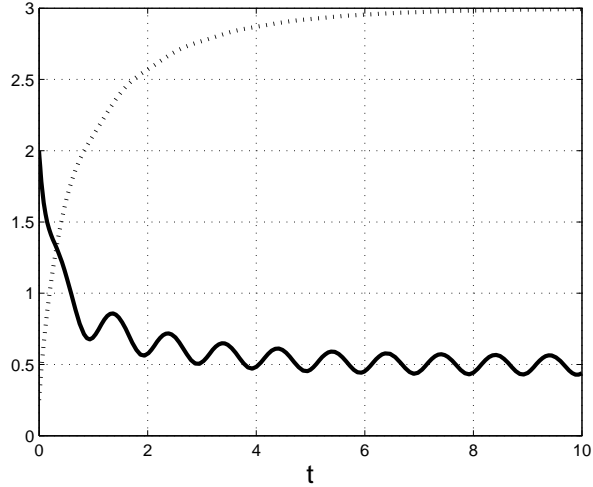}
  \caption{Solution  $x_1(t)$ (solid line)
	and~$x_2(t)$  (dashed line) of the system in Example~\ref{exa:st_ent_inf}
	as a function of~$t$.  }\label{fig:ent_inf}
  \end{center}
\end{figure}

Since~IC implies~{\stshort} and this implies~{\sostshort}, it follows from Proposition~\ref{prop:st_entrain}
that
IC   implies entrainment to~$T$-periodic vector fields.\footnote{Note that the proof that IC implies~{\stshort}
used a result for time-invariant systems, but an analogous argument holds for the time-varying case as well.}
The next example demonstrates this.
 \begin{Example}\label{exa:st_ent_inf}
Consider again the system in Example~\ref{exa:new_from_russo_inf}, and  assume that
the kinase  abundance~$u(t)$  is  a
strictly positive and periodic function of time with period~$T$.
Since we already showed that this system is~{\stshort},
it  admits a unique periodic solution~$\gamma$, of period~$T$,
and any trajectory of the system converges to~$\gamma$.
Figure~\ref{fig:ent_inf}
depicts the solution of~\eqref{eq:trans_module_inf} for~$\delta=2$, $ k_2=1$, $z_T=4$,
$e_T=3$,~$u (t)=2+\sin(2\pi t)$, and initial condition~$x_1(0)=2,x_2(0)=1/4$.
It may be seen that both state-variables converge to a periodic solution with period~$T=1$.
(In particular,~$x_2$ converges to the constant function~$x_2(t)\equiv e_T$ that is of course
periodic with period~$T$.)
\end{Example}

Contraction can be characterized using a Lyapunov-Finsler function~\cite{contra_sep}.
The next result describes  a similar characterization for~{\stshort}.
For simplicity, we state this for the time-invariant system~\eqref{eq:time_in_var_sys}.
\begin{Proposition}\label{prop:equ_lyap}
The following two  conditions are equivalent.
\begin{enumerate}[(a)]
					\item The time-invariant system~\eqref{eq:time_in_var_sys} is {\stshort} w.r.t. a  norm~$|\cdot|$.
					\label{imp_clf:a}
					\item \label{imp_clf:b} For any~$\tau>0$ there exists~$\ell=\ell(\tau)>0$ such that
					for any~$a,b \in \Omega$ and any~$c$ on the line connecting~$a$ and~$b$
					  the solution of~\eqref{eq:var_sys} with~$x(0)=c$ and~$\delta x(0)=b-a$ satisfies
					\be\label{eq:delta_cond}
										|\delta x(t+\tau)|\leq  \exp(-\ell t) |\delta x(0)|,\quad \text{for all }t\geq 0.
					\ee
\end{enumerate}
\end{Proposition}

Note that~\eqref{eq:delta_cond}
implies that the function~$V(x,\delta x):=|\delta x|$ is a \emph{generalized} Lyapunov-Finsler function in the following sense.
For any~$\tau>0$ there exists~$\ell=\ell(\tau)>0$ such that along solutions of the variational system:
\[
				V\left (x(t+\tau,x(0)),\delta x(t+\tau,\delta x(0),x(0))\right )\leq  \exp(-\ell t) V( x(0)
				,\delta x(0)),\quad \text{for all }t\geq 0.
\]

In the next section, we describe several more related notions and explore the relations between
them.

\section{Additional Notions and Relations}\label{sec:relations}

It is straightforward to show that
each of the three generalizations of  contraction in Definition~\ref{def:qcont}
  implies  that~\eqref{eq:fdyn} is~NE.
 One may perhaps expect that any of the three generalizations
 of  contraction in Definition~\ref{def:qcont}
also implies~WC.
Indeed,~{\stshort} does imply WC, because
\[
|x(s_2, s_1, a) - x(s_2, s_1, b)| \leq  \exp\left ( - \ell(s_2 - s_1)/2  \right )|a - b| < |a - b|,
\]
for all~$ 0 \leq s_1 <  s_2$ if {\stshort} holds (simply apply the definition with~$ t_1 = s_1$,
$ \tau = (s_2 - s_1)/2 > 0$, and~$t_2 = s_1 + \tau $ in~\eqref{eq:ucont}).
However, the next example shows that~{\soshort} does not imply WC.
 \begin{Example}\label{exa:eps}
 Consider the  {scalar}   system
            \be\label{eq:shift}
                    \dot{x}=\begin{cases}  -2x, & 0\leq |x| <1/2,\\
                                               -\frac{x}{|x|}, & \frac{1}{2}\leq |x|\leq 1 ,
                               \end{cases}
            \ee
            with~$x$ evolving on~$\Omega:=[-1,1]$.
Clearly, this system is not WC.
However, it is not difficult to show that it satisfies
 the definition of~{\soshort} with~$\ell =\ell(\varepsilon):=\min\{ \ln(1+\varepsilon),1\}$ .

\end{Example}


The next result presents two conditions that are equivalent to~{\sostshort}.
\begin{Lemma}\label{lem:eqdef}
The following conditions are equivalent.
 \begin{enumerate}
 \item \label{item1}
 System~\eqref{eq:fdyn}  is {\sostshort} on~$\Omega$
  w.r.t.  some vector norm $|\cdot|_v:\R^n\to \R_+  $.
\item \label{item2}
 For each~$ {\tau}>0$
  there exists~$ {\ell}= {\ell}( {\tau} )>0$
 such that
 \begin{align}\label{eq:qcontop}
            |x(t_2+\tau,  t_1,a)-  x(t_2+\tau,t_1,b)|_v
                    &\leq   (1+ {\tau}) \exp(-  (t_2-t_1)  {\ell} ) |a-b|_v ,
 \end{align}
 for all~$t_2\geq t_1\geq 0$ and all~$a,b \in \Omega$.
\item \label{item3}
For each~$\varepsilon >0$ and each~$\tau>0$
  there exists~$\ell_1=\ell_1(\tau,\varepsilon)>0$
 such that
 \be\label{eq:qcontnew}
            |x(t  ,t_1,a)-x(t  ,t_1,b)|_v  \leq   (1+\varepsilon) \exp(-  (t -t_1) \ell_1 ) |a-b|_v ,
 \ee
 for all~$t  \geq t_1+\tau \geq \tau$ and all~$a,b \in \Omega$.
\end{enumerate}
\end{Lemma}

\begin{figure*}[t]
 \begin{center}
 \input{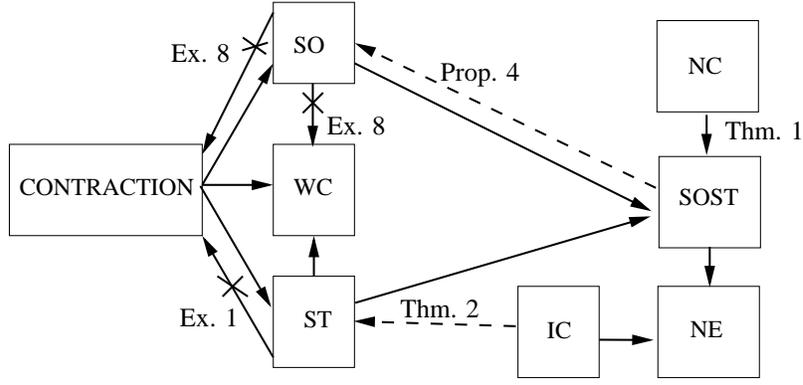}
\caption{Relations between various contraction notions.
An arrow denotes implication;
a crossed out arrow denotes that the implication is in general false;
and a dashed arrow denotes
 an implication that holds under an additional condition.
Some of the relations are immediate.
Others follow from the results   marked  near the arrows.
}
\label{fig:graphbn2}
\end{center}
\end{figure*}

Fig.~\ref{fig:graphbn2} summarizes
the relations between the various contraction
notions.

 \section{Proofs} \label{sec:proofs}
 \noindent  {\sl Proof of Theorem~\ref{thm:qcon}.}
Fix arbitrary $\varepsilon>0$ and $t_1\geq 0$.
The function $\zeta=\zeta( \tau) \in (0,1/2]$ is as in the statement of the
Theorem.
For each $ \tau>0$, let $c_{\zeta}>0$ be a contraction constant
on~$\Omega_\zeta$, where we write $\zeta=\zeta( \tau)$ here and in what follows.
Pick~$a,b \in \Omega$ and~$\tau>0$.
By~\eqref{eq:enter},
$x(t,t_1,a),x(t,t_1,b) \in \Omega_\zeta$ for  all~$t\geq t_1+\tau$, so
 \begin{align}\label{eq:adder1}
|x(t ,t_1,    a)- x(t ,t_1, b)|_\zeta
& \leq  \exp(-c_{\zeta } (t- t_1-\tau))   | x(t_1+\tau ,t_1, a)- x(t_1+\tau ,t_1,b)|_\zeta,\quad\text{for  all }t\geq t_1+\tau.
 \end{align}
In particular,
 \begin{align}\label{eq:zereg}
|x(t ,t_1,    a)- x(t ,t_1, b)|_\zeta
&< | x(t_1+\tau ,t_1, a)- x(t_1+\tau ,t_1,b)|_\zeta,\quad\text{for  all }t> t_1+\tau.
 \end{align}
From the convergence property of norms in the Theorem statement,
there exist~$v_{\zeta },w_{\zeta } >0 $ such
 that
 \begin{align}
                            |y|  & \leq v_{\zeta } | y|_{\zeta } \leq w_{\zeta } v_{\zeta } | y|,\quad \text{for all } y\in\Omega,
 \end{align}
and $v_{\zeta }\rightarrow1$, $w_{\zeta }\rightarrow1$ as $\tau\rightarrow0$.
Combining this with~\eqref{eq:zereg} yields
 \begin{align*}
|x(t ,t_1,    a)- x(t ,t_1, b)|
&< v_\zeta w_\zeta | x(t_1+\tau ,t_1, a)- x(t_1+\tau ,t_1,b)|,\quad\text{for  all }t> t_1+\tau.
 \end{align*}
Note that taking~$\tau\to   0$ yields
 \begin{align}\label{eq:eqadd}
|x(t ,t_1,    a)- x(t ,t_1, b)|
&\leq  |  a - b|,\quad\text{for  all }t> t_1 .
 \end{align}
Now for~$t\geq t_1+\tau$ let~$p:=t- t_1-\tau$. Then
  \begin{align*}
     |  x(t ,   t_1,  a) - x(t ,t_1,b)|
          & \leq  v_\zeta   |  x(t ,   t_1,  a) - x(t ,t_1,b)|_\zeta \\
          & \leq  v_\zeta  \exp(-c_\zeta p) | x(t_1+\tau ,t_1, a)- x(t_1+\tau ,t_1,b)|_\zeta \\
          &\leq v_\zeta   w_\zeta  \exp(-c_\zeta p) | x(t_1+\tau ,t_1, a)- x(t_1+\tau ,t_1,b)|\\
          &\leq v_\zeta   w_\zeta \exp(-c_\zeta p)    | a-b| ,
  \end{align*}
  where the last inequality follows from~\eqref{eq:eqadd}.
Since
$v_{\zeta }\rightarrow1$, $w_{\zeta }\rightarrow1$ as $\tau\rightarrow0$,
$v_{\zeta } w_{\zeta }  \leq  {1+ \varepsilon}$
for $\tau>0$ small enough.
Summarizing, there exists~$\tau_m=\tau_m(\varepsilon)>0$ such that for all~$\tau \in[0, \tau_m]$
  \begin{align}\label{eq:summa}
    |& x(t+\tau ,t_1,  a)- x(t+\tau ,t_1,b)|
          \leq  (1+\varepsilon)   \exp(-c_\zeta (t- t_1 ))    | a-b| ,
  \end{align}
for all~$a,b\in \Omega$ and all~$t\geq t_1 $. Now pick~$\tau > \tau_m$.
For any~$t\geq t_1$,
let~$s:= t+\tau-\tau_m $.
Then
  \begin{align*}
    |  x(t+\tau ,t_1,   a)-   x(t+\tau ,t_1,b)|
    &=|  x(s+\tau_m ,t_1,  a)- x(s+\tau_m ,t_1,b)|  \\
    &\leq  (1+\varepsilon)   \exp(-c_\zeta (s- t_1 ))    | a-b| \\
     &\leq  (1+\varepsilon)   \exp(-c_\zeta (t- t_1 ))    | a-b|,
  \end{align*}
and this  completes the proof.~$\QED$

\noindent
{\sl Analysis of the system from Example~\ref{exa:bio_smith}.}
The Jacobian  of~\eqref{eq:bio} is
\be\label{eq:jac}
            J(x)= \begin{bmatrix}     -\alpha_1 &  0       & 0 &\dots&0&  g'(x_n)\\
                                              1 & -\alpha_2& 0 &\dots& 0 &0 \\
                                              0& 1 & -\alpha_3 &\dots& 0 &0 \\
                                              &&\vdots\\
                                              0& 0 & 0  &\dots& 1 &-\alpha_n                                     
                                                 \end{bmatrix},
\ee
so
\begin{align*}
 D_\varepsilon & J(x) D_\varepsilon ^{-1} = \begin{bmatrix}
   -\alpha_1 & 0 & 0  &\dots & 0 &  \frac{ g'(x_n) }{ \prod_{i=1}^{n-1}(\alpha_{i}-\varepsilon )  } \\
  \alpha_1-\varepsilon &  -\alpha_2 &0& \dots &0 &0\\
  0 &   \alpha_2 -\varepsilon &0& \dots &0 &0\\
  &&\vdots\\
  0 &  0   &0& \dots &\alpha_{n-1}-\varepsilon  &-\alpha_n
    \end{bmatrix}.
\end{align*}
Thus,
\begin{align}\label{eq:2term}
\mu_{1,D_\varepsilon}(J(x))& =
            \max\{   -  \varepsilon, \frac{ g'(x_n)-\alpha_n \prod_{i=1}^{n-1}(\alpha_{i}-\varepsilon )}
            {\prod_{i=1}^{n-1}(\alpha_{i}-\varepsilon )  }
             \}.
\end{align}
 Suppose that~$ {k-1}  < \alpha {k^2}$. Then for all~$x\in\R^n_+$,
\[
        g'(x_n)= \frac{k-1}{(k+x_n)^2} \leq \frac{k-1}{k^2}<\alpha .
\]
Combining this with~\eqref{eq:2term} implies that
there exists a sufficiently small~$\varepsilon>0$ such that
 $\mu_{1,D_\varepsilon}(J(x)) <-\varepsilon/2$
for all~$x\in \R^n_+$,  so the system is contractive on~$\R^n_+$ w.r.t.~$|\cdot|_{1,D_\varepsilon}$.

Now assume that~${k-1} =\alpha {k^2}$.
By~\eqref{eq:jac},
\[
            \det(J(x))=(-1)^n (\alpha-g'(x_n)),
\]
so for every~$x \in \R^n_+$ with~$x_n=0$, we have~$\det(J(x))=(-1)^n (\alpha-g'(0))=0$.
This implies that the system does not satisfy~\eqref{eq:contdef}, w.r.t. any norm, on~$\R^n_+$.

We now use Theorem~\ref{thm:qcon} to
prove that~\eqref{eq:bio}  is {\sostshort}.
Since~$g'(u)=\frac{k-1}{(k+u)^2}$ and~$k>1$,
\[
            g(x_n)\geq g(0)=1/k,\quad \text{for all } x \in \R^n_+.
\]
For~$\zeta \in (0,1/2]$, let
\[
            \Omega_\zeta:=\{x \in \R^n_+: x\geq \zeta\}.
\]
It is straightforward to verify that~\eqref{eq:bio} satisfies   condition~(BR)
in~\cite[Lemma~1]{RFM_entrain}, and this implies that for every~$\tau>0$ there exists~$\varepsilon (\tau)>0$ such that
$x(t)\in \Omega_\varepsilon $ for all~$t\geq \tau$. Then
\[
        g'(x_n)= \frac{k-1}{(k+x_n)^2} \leq \frac{k-1}{(k+\varepsilon)^2}< \frac{k-1}{k^2} = \alpha .
\]
We already showed that this implies that there exists a~$\zeta>0$
and a norm~$|\cdot|_{1, D_\zeta }$
such that~\eqref{eq:bio} is contractive on~$\Omega_\varepsilon$ w.r.t. this norm.
  Summarizing, all the conditions in Theorem~\ref{thm:qcon} hold, and we conclude that
 \eqref{eq:bio}   is~{\sostshort} on~$\R_+^n$ w.r.t.~$|\cdot|_{1,D_0}$.~$\QED$

\noindent {\sl Analysis of the system in Example~\ref{exa:phos:relay}.}
For $a\in\Omega$, let~$x(t,t_1,a)$ denote the solution of~\eqref{eq:phos} at time~$t\geq t_1$
for the initial condition~$x(t_1)=a$. Pick~$\tau>0$.
Eq.~\eqref{eq:phos} satisfies   condition~(BR)
in \cite[Lemma~1]{RFM_entrain}, and this implies that
there exists~$\varepsilon=\varepsilon(\tau)>0$ such that
for all~$a\in\Omega$,    all~$i=1,\dots n$, and all~$t\geq  t_1+  \tau$
\[
x_i(t,t_1,a)\geq \varepsilon .
\]
Furthermore, if we define~$y_i(t):=p_{n-i+1}-x_{n-i+1}(t)$, $i=1,\dots,n$,
 then the~$y$ system also satisfies      condition~(BR)
in  \cite[Lemma~1]{RFM_entrain},
and this implies that   there exists~$\varepsilon_1=\varepsilon_1(\tau)>0$ such that
for all~$a\in\Omega$,    all~$i=1,\dots n$, and all~$t\geq t_1+ \tau$
\[
y_i(t,t_1,a) \geq \varepsilon_1.
\]
We conclude that after an arbitrarily short time~$\tau >0 $ every
 state-variable~$x_i(t)$, $t\geq \tau+t_1$,  is separated from~$0$ and from~$p_i$. This means
the following.
For~$\zeta \in [0,1/2]$, let
\[
            \Omega_\zeta:=\{x \in \Omega:  \; \zeta p_i \leq x_i \leq (1-\zeta) p_i,\; i=1,\dots,n\}.
\]
Note   that~$\Omega_0=\Omega$, and that~$ \Omega_\zeta$
is a strict subcube of~$\Omega$  for all~$\zeta \in (0,1/2] $.
 Then for any~$t_1 \geq 0$, and any~$\tau>0$
 there exists~$\zeta=\zeta(\tau)\in (0,1/2)$,
with~$\zeta (\tau)\to 0$ as~$\tau \to 0 $,
such that
\be\label{eq:zinomeps}
            x(t ,t_1,a) \in \Omega_\zeta,\quad \text{for all } t\geq t_1+\tau \text{ and all }a\in\Omega.
\ee

The Jacobian of~\eqref{eq:phos} satisfies~$J(t,x)=L(x)-\diag(  c(t) ,0,\dots,0,\eta_n    ) $,
where
\begin{equation*} \label{eq:floatingequation}
             L(x) = \begin{bmatrix}
  -\eta_1 (p_2-x_2 )& \eta_1  x_1                                                                       & 0        &0 \\
\eta_1 (p_2-x_2 )          &  -\eta_1  x_1 -\eta_2 (p_3-x_3 )                                      & \dots    & 0 \\
0                         &  \eta_2 (p_3-x_3 )                     &\dots & 0 \\
                          &                                   & \ddots\\
0                         & \dots                                                                 & -\eta_{n-2}  x_{n-2}  -\eta_{n-1}  (p_n-x_n )   &    \eta_{n-1}  x_{n-1}    \\
0                         & \dots                                                                 &     \eta_{n-1}  (p_n-x_n ) &        -\eta_{n-1} x_{n-1}
            \end{bmatrix}.
\end{equation*}
Note that~$L(x)$ is Metzler, tridiagonal,
 and has  zero sum columns for    all~$x\in \Omega$.
Note also that for any~$x\in\Omega_\zeta$ every entry~$L_{ij}$ on the sub- and super-diagonal of~$L$
satisfies~$ \zeta s_1  \leq L_{ij}  \leq   (1-\zeta) s_2$,
with~$s_2:=\max_i \{\eta_i  p_i\}>  s_1:=\min_i\{ \eta_i p_i\}  > 0$.

Note also that there exist~$x\in\partial \Omega$ such that~$J(x)$ is singular (e.g., when~$x_1=0$  and~$x_3=p_3$ the second column of~$J$ is all zeros), and this implies that the system does not satisfy~\eqref{eq:contdef} on~$\Omega$ w.r.t. any norm.

By~\cite[Theorem~4]{RFM_entrain}, for any~$\zeta\in(0,1/2]$ there exists
$\varepsilon=\varepsilon(\zeta)>0$, and a diagonal matrix~$D=\diag(1,q_1,q_1q_2,\dots,q_1q_2\dots q_{n-1})$,
with~$q_i=q_i(\varepsilon)>0$,
such that~\eqref{eq:phos} is contractive on~$\Omega_\zeta$
w.r.t. the  the scaled~$L_1$ norm defined by~$|z|_{1,D}=|Dz|_1$.
Furthermore, we can choose~$\varepsilon$ such that~$\varepsilon(\zeta) \to 0$ as~$\zeta\to 0$,
and~$D(\varepsilon )\to I$ as~$\varepsilon\to 0$.
 Summarizing, all the conditions in
Definition~\ref{def:NC} hold, so~\eqref{eq:phos} is NC on~$\Omega$ and applying
Theorem~\ref{thm:qcon}
  concludes the analysis.~$\QED$

\noindent {\sl Proof of Proposition~\ref{prop:new_meas_zero}.}
Without loss of generality, assume that~$S_0=\{1,\dots, k\}$,  with~$1  \leq k<n-1$,
so that~$S_-=\{k+1,\dots,n\}$. Fix~$\varepsilon \in (0,1)$.
Let~$D=\diag(d_1,\dots,d_n)$ with the~$d_i$s defined as follows.
For every~$i \in S_0$,~$d_i=1$ and~$d_{z(i)}=1-\varepsilon$. All the other~$d_i$s are   one.
Let~$\tilde J:=DJD^{-1}$. Then~$\tilde J_{ij}=\frac{d_i}{d_j}J_{ij}$.
We now calculate~$\mu_1(\tilde J)$. Fix~$j\in S_0$. Then~$d_j=1$, so
\begin{align*}
c_j(\tilde J)&=\tilde J_{jj}+\sum_{ \substack { 1\leq i \leq n\\  i \not =j}}|\tilde J_{ij}|\\
&=  J_{jj}+\sum_{ \substack { i\in S_0\\  i \not =j}}d_i| J_{ij}|+\sum_{ \substack { k\in S_-\\  k \not =j}}d_k| J_{kj}|\\
&=  J_{jj}+\sum_{ \substack { i\in S_0\\  i \not =j}}|  J_{ij}|+\sum_{ \substack {  k \in S_-\\  k  \not =j}}d_k |  J_{kj}|\\
&<c_j(J),
\end{align*}
where the inequality follows from  the fact that~$d_k\leq 1$ for all~$k$,
and   for the specific value~$k=z(j)\in S_-$ we have~$d_k=1-\varepsilon$ and~$|J_{kj}|>0$.
We conclude that for every~$j \in S_0$, $c_j(\tilde J)<c_j(J)=0$.
 It follows from property~\ref{item:sminus}) in the statement of
Proposition~\ref{prop:new_meas_zero}
   and the compactness of~$\Omega$
that there exists~$\delta>0$ such that~$c_j(J(x))<-\delta$ for all~$j \in S_-$ and all~$ x \in\Omega$,
so for~$\varepsilon>0$ sufficiently small we have~$c_j(\tilde J(x))<-\delta/2$ for all~$j \in S_-$ and all~$ x \in\Omega$.
We conclude  that for all~$\varepsilon>0$ sufficiently small,~$\mu_1(DJD^{-1})=\max_j c_j(\tilde J) <0$, i.e.
the system is contractive w.r.t.~$|\cdot|_{1,D}$.
Clearly,~$|\cdot|_{1,D} \to |\cdot|_1$ as~$\varepsilon \to 0$,
and  applying Corollary~\ref{coro:new} completes the proof.~\QED

 \noindent  {\sl Proof of Proposition~\ref{prop:sepimp}.}
Suppose that~\eqref{eq:fdyn} is  {\sostshort} w.r.t. some
norm~$|\cdot|_v$.
 Pick~$\varepsilon>0$.  Since the system is~{\sweshort}, there exists~$\tau_0=\tau_0(\varepsilon)>0$
such that
\[
                                        | x(t ,t_0,a)-x(t ,t_0,b) |_v
                                        \leq(1+ {\varepsilon}/{2})|a-b|_v, \
 \]
 for all $t \in [t_0,t_0+\tau_0]$.
Letting~$\ell_2:=\frac{1}{\tau_0}\ln(\frac {1+\varepsilon}{1+(\varepsilon/2)} )$   yields
 \be\label{eq:nnsep}
                                        | x(t ,t_0,a)-x(t ,t_0,b) |_v
                                        \leq(1+\varepsilon) \exp(-  (t -t_0) \ell_2 ) |a-b|_v,
\ee
 for all $ t \in [t_0,t_0+\tau_0]$.
By item~\ref{item3} in Lemma~\ref{lem:eqdef} there exists~$\ell_1=\ell_1(\tau_0,\varepsilon)>0$
 such that
 \[
            |x(t  ,t_0,a)-x(t  ,t_0,b)|_v  \leq   (1+\varepsilon) \exp(-  (t -t_0) \ell_1 ) |a-b|_v ,
 \]
 for all $t  \geq t_0+\tau_0$.
Combining this with~\eqref{eq:nnsep} yields
\[
            |x(t  ,t_0,a)-x(t  ,t_0,b)|_v  \leq    (1+\varepsilon) \exp(-  (t -t_0) \ell  ) |a-b|_v ,
 \]
 for all $t  \geq t_0$,
            where~$\ell:=\min\{\ell_1,\ell_2\}>0$. This proves {\soshort}.~$\QED$

\noindent  {\sl Proof  of Theorem~\ref{thm:time_invar}.}
We require the following result.
\begin{Lemma}\label{lem:strict}
If   system~\eqref{eq:time_in_var_sys} is IC
then  for each $\tau > 0$ there exists~$d=d(\tau) > 0$
 such  that
 \[
 \dist(x(t,x_0),\partial \Omega) \geq d ,\quad \text{for all } x_0 \in \Omega  \text{ and all } t \geq \tau.
 \]
\end{Lemma}

 {\sl Proof of Lemma~\ref{lem:strict}.}
 Pick  $\tau>0$ and~$ x_0 \in \Omega$.
Since~$\Omega$ is an invariant set,~$\Int( \Omega)$ is also an invariant set
 (see, e.g.,~\cite[Lemma III.6]{mcs_angeli_2003}), so~\eqref{eq:cond_a_enu} implies that
    $x(t,x_0) \not \in \partial \Omega$ for {all}~$t>0$. Since~$\partial \Omega$
		is compact,
$e_{x_0}:= \dist(x(\tau,x_0), \partial \Omega)  > 0$.
Thus, there exists a
  neighborhood~$U_{x_0}$ of $x_0$, such that
$\dist(x(\tau,y), \partial \Omega) \geq e_{x_0}/2$
for all $y \in U_{x_0}$.
Cover~$\Omega$ by such~$U_{x_0}$ sets.
 By compactness of~$\Omega$, we can pick a    finite subcover.
  Pick smallest~$e$ in this subcover, and denote
 this by~$d$. Then~$d>0$ and
we have that
$\dist(x(\tau,x_0), \partial \Omega) \geq d$ for all $x_0 \in \Omega$.
Now, pick~$t \geq \tau$.
Let $x_1: = x(t-\tau,x_0)$.  Then
\begin{align*}
\dist(x(t,x_0), \partial \Omega) &= \dist(x(\tau,x_1), \partial \Omega) \\&\geq d,
\end{align*}
and this completes the proof of Lemma~\ref{lem:strict}.~$\QED$

We can now prove  Theorem~\ref{thm:time_invar}.
We recall some definitions from the theory of convex sets.
Let~$B(x,r)$ denote the closed ball of radius~$r$ around~$x$ (in the Euclidean norm).
 Let~$K$ be a compact and convex set with~$0\in\Int(K)$.
Let~$s(K)$ denote the
\emph{inradius} of~$k$, i.e.
   the radius of the
largest ball contained in~$ K$.
For~$\lambda \in[0,s(K)]$  the \emph{inner parallel set of~$K$ at
distance~$\lambda$} is
\[
				K_{-\lambda} := \{x\in K: B(x,\lambda) \subseteq K\}.
\]
Note that~$K_{-\lambda}$ is a compact and convex set;
in fact,~$K_{-\lambda}$ is the intersection of all
the translated support hyperplanes of~$K$, with each hyperplane
translated  ``inwards'' through a distance~$\lambda$
(see~\cite[Section~17]{convex_syn}).
Assume, without loss of generality,  that~$0\in\Int(\Omega)$.
Pick~$\tau>0$.
Let~$M=M(\tau):=\{ x(t,x_0): t \geq \tau,\; x_0 \in \Omega  \}$.
By Lemma~\ref{lem:strict}, $M\subset \Omega$ and~$\dist(y,\partial \Omega)\geq d>0$ for all~$y\in M$. Let~$\lambda=\lambda(\tau):=\frac{1}{2}\min\{d,s(\Omega)\}$.
Then~$\lambda>0$.
Pick~$z \in M$. We claim that~$B(z,\lambda) \subseteq \Omega $.
To show this, assume that there exists~$v\in B(z,\lambda) $ such that~$v\not\in \Omega$.
 Then there is a point~$q$ on the line connecting~$v$ and~$z$
such that~$q\in \partial \Omega$. Therefore,
\begin{align*}
											\dist(z,\partial \Omega)&\leq |z-q|\\
																				 &\leq |z-v|\\
																				 &\leq \lambda\\
																				 &\leq d/2,
\end{align*}
and this is a contradiction as~$z\in M$.
We conclude that~$M\subseteq K_{-\lambda}$.
 Let~$c=c(\tau):=\max_{x \in  K_{-\lambda}} \mu(J(x))$. Then~\eqref{eq:mucom_inv} implies that~$c <0$. Thus, the system
 is contractive on~$ K_{-\lambda} $, and
for all~$a,b \in \Omega$ and all~$t\geq 0$
\[
            |x(t+\tau,a)-x(t+\tau,b)|\leq \exp( c  t) |a-b|,
\]
 where~$|\cdot| $
is the vector norm corresponding to the matrix measure~$\mu $.
  This establishes~\stshort, and thus
  completes the proof of Theorem~\ref{thm:time_invar}.~$\QED$


\noindent  {\sl Proof of Corollary~\ref{coro:attract}.}
\noindent Since~$\Omega$ is convex, compact, and invariant, it includes an
 equilibrium point~$e$  of~\eqref{eq:time_in_var_sys}. Clearly,~$e\in\Int(\Omega)$.
By Theorem~\ref{thm:time_invar},
 the system is~{\stshort}. Pick~$a \in \Omega$ and~$\tau>0$, and let~$\ell=\ell(\tau)>0$.
 Applying~\eqref{eq:ucont}
with~$b=e$ yields
\begin{align*}
            |x(t +\tau, a)-e|
            \leq     \exp(-   \ell t ) |a-e|,
 \end{align*}
 for all~$t\geq  0$. Taking~$t \to\infty$
  completes the proof.~$\QED$
\begin{Remark}
Another possible proof of  Corollary~\ref{coro:attract}  is based
on defining~$V:\Omega \to\R_+$ by~$V(x):=|x-e|$.  Then for any~$a\in\Omega$,~$V(x(t,a))$ is nondecreasing, and the LaSalle invariance principle tells us that~$x(t,a)$ converges to an invariant
subset of the set~$\{y\in \Omega:|y-e|=r\}$, for some~$r\geq 0$.
If~$r=0$ then we are done. Otherwise, pick~$y$ in the omega limit set of the trajectory.
Then~$y\not \in \partial \Omega$, so~\eqref{eq:mucom_inv} implies that~$V$ is strictly
decreasing. This contradiction completes the proof.
\end{Remark}

\noindent
{\sl Proof of Proposition~\ref{prop:equ_lyap}.}
Pick~$a,b\in\Omega$. Let~$\gamma:[0,1]\to \Omega$ be the line~$\gamma(r):=(1-r) a+r b$.
Note that since~$\Omega$ is convex,~$\gamma(r)\in\Omega$ for all~$r\in[0,1]$.
Let
\[
			w(t,r):= \frac{d}{dr}x(t,\gamma(r)).
\]
This  measures the sensitivity of the solution at time~$t$ to a   change in the
 initial condition along the line~$\gamma$. Note that~$w(0,r)=\frac{d}{dr}\gamma(r)=b-a$, and
\[
				\dot w(t,r)=J  (x(t,\gamma(r))) w(t,r).
\]
Comparing this to~\eqref{eq:var_sys}  implies
 that~$w(t,r) $ is equal to the second component,~$\delta x(t)$, of the
solution  of the variational system~\eqref{eq:var_sys}
with initial condition
\begin{align}\label{eq:i nit_var}
x(0)&=(1-r) a+r b,\\
\delta x(0)&=b-a.\nonumber
\end{align}

Suppose that the time-invariant system~\eqref{eq:time_in_var_sys} is \stshort. Pick~$\tau>0$.
Let~$\ell=\ell(\tau)>0$. Then for any~$r\in[0,1)$ and any~$\varepsilon \in[0,1-r]$,
 \[
            |x(t +\tau,  \gamma(r+\varepsilon))-x(t +\tau ,\gamma(r ) )|  \nonumber
            \leq     \exp(- t \ell ) |\gamma(r+\varepsilon)-\gamma(r)|.
 \]
Dividing both sides of this inequality by~$\varepsilon$ and taking~$\varepsilon \downarrow  0$ implies
 that
 \be\label{eq:condwexp}
            |w(t +\tau,  r) |
            \leq     \exp(-  t \ell ) |b-a|,
 \ee
so
 \[
            |\delta x(t +\tau) |
            \leq     \exp(-  t \ell ) |\delta x(0)|.
 \]
This proves the implication~\eqref{imp_clf:a} $\to$ \eqref{imp_clf:b}.
To prove the converse implication, assume that~\eqref{eq:delta_cond} holds. Then~\eqref{eq:condwexp} holds and thus
\begin{align*}
									 |x(t+\tau,b)-x(t+\tau,a)|&=\left| \int_0^1 \frac{d}{dr} x(t+\tau,\gamma(r)) \dd r \right |\\
									&\leq  \int_0^1\left  |w(t+\tau, r ) \right |\dd r\\
									&\leq  \int_0^1 \exp(-\ell t) |b-a|\dd r\\
									& =    \exp(-\ell t) |b-a| ,
\end{align*}
so the system is~{\stshort}.~\QED

\noindent
 {\sl Proof of Lemma~\ref{lem:eqdef}.}
If~\eqref{eq:fdyn}  is {\sostshort} then~\eqref{eq:qcontop} holds  for the particular case~$\varepsilon=\tau$
in Definition~\ref{def:qcont}. To prove the converse
implication, assume that~\eqref{eq:qcontop} holds.
Pick~$\hat{\tau},\hat{\varepsilon}>0$. Let \be \label{eq:deftaumin} \tau:=\min\{\hat{\tau},\hat{\varepsilon}    \}, \ee and
let~$\ell=\ell(\tau)>0$. Pick~$t\geq t_1\geq 0$, and let~$t_2:=t+\hat{\tau}-\tau \geq t_1$. Then
 \begin{align*}
            |x(t_2 + {\tau},  t_1,a)-  x(t_2+\tau  ,t_1,b)|_v
                    &\leq   (1+  {\tau}) \exp(-  ( t_2   -  t_1)  {\ell} ) |a-b|_v \\
                    &\leq (1+\hat{\varepsilon}) \exp(-  ( t   -  t_1)  {\ell} ) |a-b|_v,
 \end{align*}
where the last inequality follows from~\eqref{eq:deftaumin}. Thus,
 \begin{align*}
            |x(t +\hat{\tau}, t_1,a)-  x(t +\hat{\tau},t_1,b)|_v
                    &\leq (1+\hat{\varepsilon}) \exp(-  ( t   -  t_1)  {\ell} ) |a-b|_v,
 \end{align*}
  and recalling that~$\hat{\tau},\hat{\varepsilon}>0$ were arbitrary,
 we conclude that Condition~2) in Lemma~\ref{lem:eqdef} implies
 \sostshort.

To prove that Condition~3) is equivalent to {\sostshort}, suppose that~\eqref{eq:qcontnew} holds. Then for any~$t_2\geq t_1$,
 \begin{align*}
                |x(t_2+\tau,t_1,a)- x(t_2+\tau,t_1,b)|_v &  \leq (1+\varepsilon) \exp(-  (t_2+\tau -t_1) \ell_1 ) |a-b|_v \\
                                                        &\leq (1+\varepsilon) \exp(-  (t_2  -t_1) \ell_1 ) |a-b| _v,
\end{align*}
so we have {\sostshort}.
 Conversely, suppose that~\eqref{eq:fdyn}  is {\sostshort}.
 Pick any~$\tau,\varepsilon>0$. Then
 there exists~$\ell=\ell( {\tau}, {\varepsilon}/2  )>0$
 such that for any~$t \geq t_1+{\tau}$
 \begin{align*}
              |x(t ,t_1,a)- x(t ,t_1,b)| _v
                                       & =  |x(t-{\tau}+{\tau},t_1,a)-  x(t-{\tau}+{\tau}  ,t_1,b)| _v  \\
              & \leq (1+ {\varepsilon}/2) \exp(-  (t-{\tau}  -t_1) \ell  ) |a-b|_v.
 \end{align*}
 Thus, for any~$c \in(0,1)$
 \begin{align*}
              |x(t ,t_1,a)-  x(t ,t_1,b)| _v
              & \leq (1+ {\varepsilon}/2) \exp( \tau c \ell ) \exp(-  (t   -t_1) c \ell  ) |a-b|_v.
 \end{align*}
Taking~$c>0$ sufficiently
 small such that~$(1+ {\varepsilon}/2) \exp( \tau c \ell ) \leq 1+\varepsilon$ implies that~\eqref{eq:qcontnew} holds for~$\ell_1:=c \ell$.
  This completes the proof
  that~\eqref{eq:qcontnew}
  is equivalent to {\sostshort}.~$\QED$

\subsubsection*{Acknowledgments}
We thank Zvi Artstein and George Weiss for   helpful comments. \textcolor{black}{We are grateful to the anonymous reviewers
 for many helpful comments}.

\end{document}